\documentclass[graybox]{svmult}

\usepackage{mathptmx}       
\usepackage{helvet}         
\usepackage{courier}        
\usepackage{type1cm}        
\usepackage{makeidx}         
\usepackage{graphicx}        
\usepackage{multicol}        
\usepackage[bottom]{footmisc}

\makeindex             

\usepackage{benstyle_ICOS}
\usepackage{amssymb,amsmath,amsfonts,color,multirow}

\newcommand*{\vspa}{\vspace{-0.0pc}}

\begin{document}

\title*{Recovering piecewise smooth functions from nonuniform Fourier measurements}
\titlerunning{Recovery from nonuniform Fourier measurements}
\author{Ben Adcock \and Milana Gataric \and Anders C. Hansen}
\institute{Ben Adcock \at Department of Mathematics, Simon Fraser University, BC V5A 1S6, Canada \email{ben\_adcock@sfu.ca}
\and Milana Gataric \at CCA, Centre for Mathematical Sciences, University of Cambridge, CB3 0WA, UK \email{m.gataric@maths.cam.ac.uk}
\and Anders C. Hansen \at DAMTP, Centre for Mathematical Sciences, University of Cambridge, CB3 0WA, UK \email{ach70@cam.ac.uk}}
%
%
\maketitle

\abstract*{}

\abstract{In this paper, we consider the problem of reconstructing piecewise smooth functions to high accuracy from nonuniform samples of their Fourier transform.  We use the framework of nonuniform generalized sampling (NUGS) to do this, and to ensure high accuracy we employ reconstruction spaces consisting of splines or (piecewise) polynomials.  We analyze the relation between the dimension of the reconstruction space and the bandwidth of the nonuniform samples, and show that it is linear for splines and piecewise polynomials of fixed degree, and quadratic for piecewise polynomials of varying degree.
}

\vspa
\section{Introduction}
\label{sec:intro}
In a number of applications, including Magnetic Resonance Imaging (MRI), electron microscopy and Synthetic Aperture Radar (SAR), measurements are collected nonuniformly in the Fourier domain.  The corresponding sampling patterns may be highly irregular; for example, one may sample more densely at low frequencies and more sparsely in high frequency regimes.  Standard tools for reconstruction from such data such as gridding \cite{JacksonEtAlGridding} seek to compute approximations to the harmonic Fourier modes, which can be then further postprocessed by conventional filtering and/or edge detection algorithms.  However, gridding methods are low order, and lead to both physical (e.g.\ Gibbs phenomena) and unphysical artefacts \cite{GelbNonuniformFourier}.

In this paper we consider high-order, artefact-free methods for the reconstruction of one-dimensional piecewise smooth functions.  To do this, we use the recently-introduced tool of nonuniform generalized sampling (NUGS) \cite{BAMGACHNonuniform1D}.  NUGS is reconstruction framework for arbitrary nonuniform samples which allows one to tailor the reconstruction space to suit the function to be approximated.  Critically, in NUGS the dimension of the reconstruction space, which we denote by $\rT$, is allowed to vary in relation to the \textit{bandwidth} $K$ of the samples.  By doing so, one obtains a reconstruction which is numerically stable and quasi-optimal.  Hence, if $\rT$ is chosen appropriately for the given function -- for example, a polynomial or spline space for smooth functions, or a piecewise polynomial space for piecewise smooth functions -- one obtains a rapidly-convergent approximation.  

The key issue prior to implementation is to determine such scaling.  In principle, this depends on both the nature of the nonuniform samples \textit{and} the choice of reconstruction space.  In this paper we provide a general analysis which allows one to simultaneously determine such scaling for all possible nonuniform sampling schemes by scrutinizing two intrinsic quantities $\zeta$ and $\gamma$ of the reconstruction space $\rT$, related to the maximal uniform growth of functions in $\rT$ and the maximal growth of derivatives in $\rT$ respectively.  Provided these are known (as is the case for many choices of $\rT$), one can immediately estimate this scaling.  As a particular consequence, for trigonometric polynomials, splines and piecewise algebraic polynomials (with fixed polynomial degree), we can show that this scaling is linear, and for piecewise algebraic polynomials with varying degree we show that it is quadratic.  The asymptotic order of such estimates is provably optimal.

\vspa 
\section{Nonuniform generalized sampling}
\label{sec:setup}
Throughout we work in the space $\rH = \rL^2(0,1)$ with its usual inner product $\ip{\cdot}{\cdot}$ and norm $\nm{\cdot}$.  Define the Fourier transform by
$\hat{f}(\omega) = \int^{1}_{0} f(x) \E^{-2\pi \I \omega x} \D x$ for $\omega \in \mathbb{R}$.
We let $\{ \Omega_N  \}_{N \in \bbN}$ be a sequence of ordered nonuniform sampling points, i.e.\ $\Omega_N =  \{ \omega_{n,N} \}^{N}_{n=1} \subseteq \bbR$ where $-\infty < \omega_{1,N} < \omega_{2,N} < \ldots < \omega_{N,N} < \infty$, and let $\{ \rT_M \}_{M \in \bbN}$ be a sequence of finite-dimensional subspaces of $\rH$.  We make the natural assumption that the sequence of orthogonal projections $\cP_{M} = \cP_{\rT_M} : \rH \rightarrow \rT_M$ converge strongly to the identity operator $\cI$ on $\rH$.  That is, any function $f \in \rH$ can be approximated to arbitrary accuracy from $\rT_M$ for sufficiently large $M$.

Our goal is the following: given the samples $\{ \hat f(\omega_{n,N}) \}^{N}_{n=1}$ compute an approximation $f_{N,M}$ to $f$ from the subspace $\rT_M$.  Proceeding as in \cite{BAMGACHNonuniform1D}, we do this via the following weighted least-squares:
\be{
\label{weightedFS}
f_{N,M}= \underset{g \in \mathrm{T_M}}{\operatorname{argmin}} \sum^{N}_{n=1} \mu_{n,N} \left | \hat{f}(\omega_{n,N} ) - \hat{g}(\omega_{n,N} ) \right |^2,
}
where $\mu_{n,N} \geq 0$ are appropriate weights (see later).  As discussed in \cite{BAMGACHNonuniform1D}, the key is to choose $M$ suitably small for a given $N$ (or equivalently $N$ suitably large for a given $M$) so that the  approximation $\{ \hat f(\omega_{n,N}) \}^{N}_{n=1} \mapsto f_{N,M} \in \rT_M$ is numerically stable and quasi-optimal.  To this end, the following estimates were shown in \cite{BAMGACHNonuniform1D}:
\begin{equation}
\label{GSbounds}
\|f-f_{N,M}\| \leq C(N,M) \inf_{g\in\mathrm{T}_M} \|f-g\|, \quad   \|f_{N,M}\| \leq C(N,M) \|f\|, \quad \forall f\in \rH,
\end{equation}
where $C(N,M)=\sqrt{C_2(N,M) / C_1(N)}$  and $C_1(N,M)$ and $C_2(N)$ are the optimal constants in the inequalities
\begin{align*}
\sum^{N}_{n=1} \mu_{n,N} \left | \hat{f}(\omega_{n,N}) \right|^2 &\geq C_1(N,M) \|f\|^2, \quad \forall f\in\mathrm{T}_M,\\
\sum^{N}_{n=1} \mu_{n,N} \left | \hat{f}(\omega_{n,N}) \right|^2 &\leq C_2(N) \|f\|^2, \quad \forall f\in\mathrm{H}.
\end{align*}
In particular, $f_{N,M}$ exists uniquely for any $f\in\mathrm{H}$ if and only if $C_1(N,M)>0$.  

\rem{
Recently, a number of other works have investigated the problem of high-order reconstructions from nonuniform Fourier data.  In \cite{GelbHines2011Frames,GelbNonuniformFourier} spectral reprojection techniques were used for this task, and a frame-theoretic approach was introduced in \cite{GelbSongFrameNFFT}.  Recovering the Fourier transform to high accuracy was studied in \cite{GelbPlatteEdge}, and in \cite{GelbHinesEdge,MartinezEtAlEdge} the problem of high-order edge detection was addressed.  A more detailed discussion is beyond the scope of this paper.  However, we note that the methods we consider in this paper based on NUGS can be shown to achieve optimal convergence rates amongst all stable, classically convergent algorithms \cite{BAACHOptimality,AdcockHansenShadrinStabilityFourier}.
}

\vspa 
\section{A sufficient condition for stability and quasi-optimality}
To ensure that $C(N,M)$ is small and finite, and hence guarantee stability and quasi-optimality via \R{GSbounds}, we first need the following density assumption:

\defn{
The sequence $\{ \Omega_N \}_{N \in \bbN}$ is uniformly $\delta$-dense for some $0 < \delta < 1$ if: (i) there exists a sequence $\{ K_N \}_N \subseteq  [0,\infty)$ with $K_N \rightarrow \infty$ as $N \rightarrow \infty$ such that $\Omega_N \subseteq [-K_N,K_N]$, and (ii) for each $N$, the density condition
$\max_{n=0,\ldots,N} \{ \omega_{n+1,N} - \omega_{n,N} \} \leq \delta$
holds, where $\omega_{0,N} = \omega_{N,N} - 2 K_N$ and $\omega_{N+1,N} = \omega_{1,N} + 2 K_N$.
}
This condition ensures that the sample points spread to fill the whole real line whilst remaining sufficiently dense.\footnote{We remark in passing that the case of critical density $\delta = 1$ can also be addressed \cite{BAMGACHNonuniform1D}, but one cannot in general expect stable reconstruction for $\delta > 1$.  See also \cite{GrochenigIrregular,GrochenigIrregularExpType}.}  We will commonly refer to the numbers $K_N$ as the sampling \textit{bandwidths}.  Note that the $\delta$-dense sample points can have arbitrary locations.  In particular, the points $\{ \omega_{n,N} \}^{N}_{n=1}$ are allowed to cluster arbitrarily.  To compensate for this, we choose the weights $\mu_{n,N}$ in the least-squares \R{GSbounds} as follows:
\be{
\label{weights}
\mu_{n,N} = \tfrac12 \left(\omega_{n+1,N} - \omega_{n-1,N} \right ),\quad n=1,\ldots,N.
}
With this to hand, we next define the \textit{$z$-residual} of a finite-dimensional $\rT \subseteq \rH$:
\begin{equation*}
E_{\rT}(M,z) = \sup \left \{ \| \hat{f} \|_{\mathbb{R} \backslash (-z,z)} : \ f \in \mathrm{T}_M, \| f \| =1  \right \}, \quad z\in (0,\infty).
\end{equation*}
Here $\nm{f}_{I} = \sqrt{\int_{I} | f(x) |^2 \D x}$ denotes the Euclidean norm over a set $I$.

\begin{theorem}[\cite{BAMGACHNonuniform1D}]
\label{t:Cbound_arbitrary}
Let $\{ \Omega_N \}_{N \in \bbN}$ be uniformly $\delta$-dense, $\{ \rT_M \}_{M \in \bbN}$ be a sequence of finite-dimensional subspaces and let $0 < \epsilon < 1-\delta$.  Let $M,N \in \bbN$ be such that
\be{
\label{E_condition}
E_{\rT}(M,K_N-1/2)^2 \leq \epsilon(2-\epsilon),
}
then the reconstruction $f \mapsto f_{N,M}$ defined by \R{weightedFS} with weights given by \R{weights} has constant $C(N,M)$ satisfying
\be{
\label{Cs_bound}
C(N,M) \leq \frac{1+\delta}{1-\epsilon - \delta}.
}
\end{theorem}
This theorem reinterprets the required scaling of $M$ and $N$ in terms of the $z$-residual $E(M,K_N-1/2)$.  Note that this residual is independent of the geometry of the sampling points, and depends solely on bandwidths $K_N$.  Hence, provided \R{E_condition} holds, one ensures stable, quasi-optimal recovery for \textit{any} sequence of sample points $\{ \Omega_N \}_{N \in \bbN}$ with the same parameters $K_N$.

Unsurprisingly, the behaviour of the $z$-residual depends completely on the choice of subspaces $\{ \rT_M \}_{M \in \bbN}$.  Whilst one can often derive estimates for this quantity using ad-hoc approaches for each particular choice of $\{ \rT_M \}_{M \in \bbN}$ -- for example, see \cite{BAMGACHNonuniform1D,AHPWavelet} for the case of wavelet spaces -- it is useful to have a more unified technique to reduce the mathematical burden.  We now present such an approach.

\defn{[\cite{SzyldOblProj}]
Let $\rU$ and $\rV$ be closed subspaces of $\rH$ with corresponding orthogonal projections $\cP_{\rU}$ and $\cP_{\rV}$ respectively.  The gap between $\rU$ and $\rV$ is the quantity $G(\rU,\rV) = \| (\cI - \cP_{\rU} ) \cP_{\rV} \|$, where $\cI : \rH \rightarrow \rH$ is the identity.
}

\lem{
\label{z-residual_gap}
Let $\{ \rT_M \}_{M \in \bbN}$ and $\{ \rS_L \}_{L \in \bbN}$ be sequences of finite-dimensional subspaces of $\rH$.  Then $E_{\rT}(M,z) \leq E_{\rS}(L,z)+ G(\rS_L,\rT_M)$ for every $M,L \in \bbN$. 
}
\prf{
Let $f \in \rT_M$, $\| f \| = 1$.  Then
\eas{
\| \hat{f} \|_{\bbR \backslash (-z,z)} 
& \leq  \| \widehat{\cP_{\rS_L} f} \|_{\bbR \backslash (-z,z)} + \| f - \cP_{\rS_L} f \|
\\
& \leq E_{\rS}(L,z) \| \cP_{\rS_L} f \| + G(\rS_L,\rT_M) \| f \| \leq E_{\rS}(L,z) + G(\rS_L,\rT_M).
}
}
This lemma implies the following: if the behaviour of $z$-residual $E_{\rS}(L,z)$ and the gap $G(\rS_L,\rT_M)$ are known, then one can immediately determine the required scaling of $M$ with $z$ to ensure that $E_{\rT}(M,z)$ satisfies \R{E_condition}.   We now make the following choice for $\{ \rS_L \}_{L \in \bbN}$ to allow us to exploit this lemma:
\be{
\label{pcwse_const}
\rS_{L} = \left \{ g \in \rH : g |_{[l/L,(l+1)/L)} \in \bbP_0,\ l=0,\ldots,L-1 \right \}.
}
Here $\bbP_0$ is space of polynomials of degree zero.  In \cite{BAMGACHNonuniform1D}, it was shown that there exists a constant $c_0(\epsilon) > 0$ such that $E_{\rS}(L,z) \leq \epsilon$ whenever $z \geq c_{0}(\epsilon) L$.  Therefore, according to Lemma \ref{z-residual_gap}, to estimate $E_{\rT}(M,z)$ we now only need to determine $G(\rS_L,\rT_M)$.  

From now on, we let $0  < w_1 < \ldots < w_k < 1$ be a fixed sequence of nodes, and define the space $\rH^1_w(0,1) =  \{ f : f |_{(w_j,w_{j+1})} \in \rH^1(w_j,w_{j+1}),\  j=0,\ldots,k \}$ where $w_0 = 0$, $w_{k+1} = 1$ and $\rH^1(I)$ is the usual Sobolev space of functions on an interval $I$.  By convention, if $k=0$ then $\rH^1_w(0,1) = \rH^1(0,1)$. 
\lem{
\label{l:Gap}
Suppose that $\rT_M \subseteq \rH^1_w(0,1)$ and let $\rS_{L}$ be given by \R{pcwse_const}.  If $L^{-1} \leq \eta = \min_{j=0,\ldots,k} \{ w_{j+1} - w_j \}$ then $G(\rS_L,\rT_M) \leq \sqrt{ \gamma^2_M/(\pi L)^2 + 4 \zeta^2_M/L}$, where
\eas{
\gamma_{M} &= \max_{j=0,\ldots,k} \sup \left \{ \| f' \|_{(w_j,w_{j+1})} : f \in \rT_M, \| f \|_{(w_j,w_{j+1})} =1 \right \},
\\
\zeta_M &= \max_{j=0,\ldots,k}  \sup \left \{ \| f \|_{\infty,(w_j,w_{j+1})} : f \in \rT_M, \| f \|_{(w_j,w_{j+1})} =1 \right \},
}
and, if $I$ is an interval, $\| f \|^2_{I} = \int_{I} |f(x) |^2 \D x$ and $\| f \|_{\infty,I} = \mathrm{ess}\,\mathrm{sup}_{x \in I} |f(x) |$.  Moreover, if $k=0$, i.e.\ $\rT_M \subseteq \rH^1(0,1)$, then  $G(\rS_L,\rT_M) \leq \gamma_M / (\pi L)$.
}
\prf{
Since $L \geq 1/\eta$ there exist $l_{j} \in \bbN$ with $l_1 < l_2 < \ldots < l_{k}$ such that $0 \leq L w_{j} - l_j < 1$ for $j=1,\ldots,k$.  For an interval $I \subseteq \bbR$, let us now write $f_{I} = \frac{1}{\abs{I}} \int_{I} f$.  Then
\eas{
\| f- \cP_{\rS_L} f \|^2 =& \sum^{L-1}_{l=0} \int_{I_l} \left | f - f_{I_l} \right |^2 = \sum^{L-1}_{\substack{l=0 \\ l \neq l_1,\ldots,l_k}}  \int_{I_l} \left | f - f_{I_l} \right |^2 + \sum^{k}_{j=1} \int_{I_{l_j}} \left | f - f_{I_{l_j}} \right |^2,
}
where $I_l = [l/L,(l+1)/L)$.  Since $f \in \rH^1(I_l)$ for $l \neq l_1,\ldots,l_k$, an application of Poincar\'e's inequality gives that
\be{
\label{hudson}
\| f - \cP_{\rS_L} f \|^2 \leq \frac{1}{(L \pi )^2} \sum^{L-1}_{\substack{l=0 \\ l \neq l_1,\ldots,l_k}} \| f' \|^2_{I_l} + \sum^{k}_{j=1} \int_{I_{l_j}} \left | f - f_{I_{l_j}} \right |^2.
}
We now consider the second term.  Write
$I_{l_j} = (l_j/L,w_j) \cup (w_j,(l_j+1)/L) =  A_j \cup B_j$ and note that for an arbitrary interval $I$ we have $\int_{I} \left | f - f_I \right |^2 = \| f \|^2_{I} - | I | |f_I|^2$.  Hence
\eas{
\int_{I_{l_j}} \left | f - f_{I_{l_j}} \right |^2 & = \int_{A_j} \left | f - f_{A_j} \right |^2 + \int_{B_j} \left | f - f_{B_j} \right |^2 + \frac{|A_j| | B_j | }{|A_j|+|B_j|} \left | f_{A_j}- f_{B_j} \right |^2
\\
& \leq \frac{1}{(\pi L)^2} \left ( \| f' \|^2_{A_j} + \| f' \|^2_{B_j} \right ) + \frac{2 |A_j| | B_j | }{|A_j|+|B_j|} \left ( \| f \|^2_{\infty,A_j} + \| f \|^2_{\infty,B_j} \right ),
}
where in the final step we use Poincar\'e's inequality once more and the fact that $f$ is $\rH^1$ within $A_j$ and $B_j$.  Since $|A_j|,|B_j| \leq L^{-1}$ and $|A_j| + |B_j| = | I_{l_j} | = L^{-1}$ we now get
\eas{
 \sum^{k}_{j=1} & \int_{I_{l_j}} \left | f - f_{I_{l_j}} \right |^2  \leq \frac{1}{(\pi L)^2} \sum^{k}_{j=1} \left ( \| f' \|^2_{A_j} + \| f' \|^2_{B_j} \right ) + \frac{4}{L} \sum^{k}_{j=0} \| f \|^2_{\infty,(w_j,w_{j+1}) }.
}
Combining this with \R{hudson} now gives that
\eas{
\| f - \cP_{\rS_L} f \|^2 \leq \left ( \frac{\gamma_M}{L \pi } \right )^2 \sum^{k}_{j=0} \| f \|^2_{(w_j,w_{j+1})} + \frac{4 \zeta^2_M}{L} \sum^{k}_{j=0} \| f \|^2_{(w_j,w_{j+1})}.
}
Since $\| f \|^2 = \sum^{k}_{j=0}  \| f \|^2_{(w_j,w_{j+1})}$ the result now follows.
}
This lemma provides the main result of this paper.  Using it, we deduce that for any $\{ \rT_M \}_{M \in \bbN}$, the question of stable reconstruction now depends solely on the quantities $\gamma_M$ and $\zeta_M$, which are intrinsic properties of the subspaces completely unrelated to the sampling of the Fourier transform.

\vspa 
\section{Examples}
To illustrate this result, we end by presenting several examples.

\vspace{0.5pc}\noindent
\textbf{Trigonometric polynomials.}  Functions $f$ that are smooth and periodic can be approximated in finite-dimensional spaces of trigonometric polynomials $\mathrm{T}_M = \left \{ \sum^{M}_{m=-M} a_m \E^{2\pi\I m x}  : a_m \in \mathbb{C} \right \}$.  If $f \in \mathrm{C}^{\infty}(\mathbb{T})$, where $\mathbb{T} = [0,1)$ is the unit torus, then the projection error $\| f - \mathcal{P}_{\mathrm{T}_M f} \|$ decay superalgebraically fast in $M$; that is, faster than any power of $M^{-1}$.  If $f$ is also analytic then the error decays exponentially fast.

For this space, we have $\rT_M \subseteq \rH^1(0,1)$ and $\gamma_M \leq 2 \pi M$ by Bernstein's inequality.  Hence Theorem \ref{t:Cbound_arbitrary} and Lemmas \ref{z-residual_gap} and \ref{l:Gap} give that the reconstruction $f_{N,M}$ is stable and quasi-optimal provided $M$ scales linearly with the sampling bandwidth $K_N$.  This result extends a previous result of \cite{BAACHOptimality} to the case of arbitrary nonuniform samples.  Note that this is the best scaling possible up to a constant: for an arbitrary sequence $\{ \rT_M \}_{M \in \bbN}$ with $\dim(\rT_M) = M$ the scaling of $M$ with $K_N$ is at best linear \cite{BAMGACHNonuniform1D}.

\vspace{0.5pc} \noindent
\textbf{Algebraic polynomials.}  Functions that are smooth but nonperiodic can be approximated by algebraic polynomials.  If $\rT_M = \bbP_{M}$ is the space of algebraic polynomials of degree at most $M$, then the projection error $\| f - \mathcal{P}_{\mathrm{T}_M f} \|$ decays superalgebraically fast in $M$ whenever $f \in \rC^{\infty}[0,1]$, and exponentially fast when $f$ is analytic.

The classical Markov inequality for this space gives that $\gamma_M \leq \sqrt{2} M^2$, $\forall M \in \bbN$ \cite{BoettcherMarkov}.  Hence we deduce stability and quasi-optimality of the reconstruction, but only with the square-root scaling $M = \ord{\sqrt{K_N}}$, $N \rightarrow \infty$ (this result extends previous results \cite{BAACHAccRecov,AdcockHansenSpecData,hrycakIPRM} to the case of nonuniform Fourier samples).  On the face of it, this scaling is unfortunate since it means the approximation accuracy of $f_{N,M}$ is limited to root-exponential in $K_N$, which is much slower than the exponential decay rate of the projection error.  However, such scaling is the best possible: as shown in \cite{AdcockHansenShadrinStabilityFourier}, any reconstruction algorithm (linear or nonlinear) that achieves faster than root-exponential accuracy for analytic functions must necessarily be unstable.

\begin{figure}[t] 
\vspace{-0.2cm}
\footnotesize
\begin{center}
$\begin{array}{cccc}
\ \includegraphics[scale=0.58]{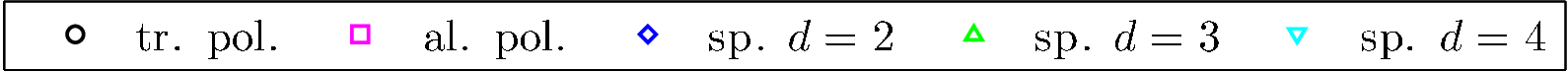} \\
\includegraphics[scale=0.655]{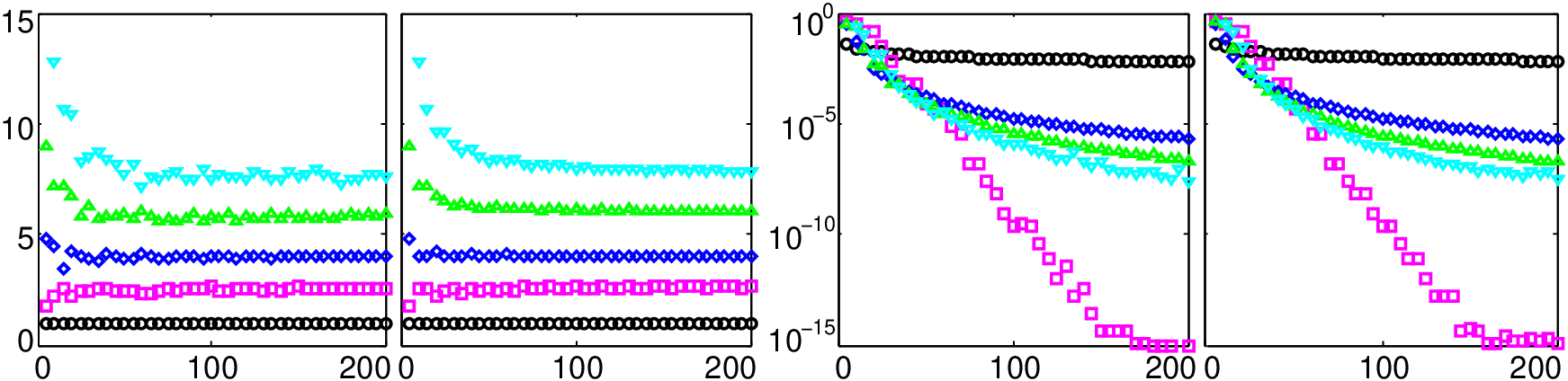} 
\end{array}$
\vspace{-0.2cm}
\caption{In the first pair of panels, depending on the type of the reconstruction space, appropriate ratios are shown: $M/K_N$ (for trigonometric polynomials), $M/\sqrt{K_N}$ (for algebraic polynomials) and $Md^2/K_N$ (for splines of order $d$), where for a given $K_N\in[5,200]$,  we used $M = \max\{ M \in \bbN: C(N,M)\leq 3\}$. In the second pair of panels, for such $K_N$ and $M$, the error $\|f-f_{N,M}\|$ is plotted  where $f(x)=x^2+x\sin(4\pi x)-\exp(x/2)\cos(3\pi x)^2$. We used different sampling schemes $\Omega_N$: jittered (for the first and third panel) and log (for the second and forth panel).}
\label{fig:rate_and_error_s}
\end{center}
\vspace{-0.7cm}
\end{figure}

\vspace{0.5pc} \noindent
\textbf{Piecewise algebraic polynomials.}  There are two issues with the previous result.  First,  the space is not suitable for approximating piecewise smooth functions.    Second, the scaling is severe.  To mitigate both issues, we may consider spaces of piecewise polynomials on subintervals.  In the first case, we fix the intervals corresponding to the discontinuities of the function, and vary the polynomial degree.  In the second case, we vary the subinterval size whilst keeping the polynomial degree fixed.

Mathematically, both scenarios equate to considering the subspaces $\rT_{w,M} = \{ f \in \mathrm{H} : f |_{[w_{j},w_{j+1})} \in \mathbb{P}_{M_j},\  j=0,\ldots,k \}$, where $w=\{ w_1,\ldots,w_{k} \}$ for $0 = w_0 < w_1 < \ldots w_k < w_{k+1} = 1$ and $M = \{ M_0,\ldots,M_{k} \} \in \mathbb{N}^{k+1}$.  If $f$ is piecewise smooth with jump discontinuities at known locations $0 = w_0 < w_1 < \ldots w_k < w_{k+1} = 1$ then the projection error decays superalgebraically fast in powers of $(M_{\min})^{-1}$ as $M_{\min}$ increases, where $M_{\min} = \min \{ M_0,\ldots,M_k \}$, and exponentially fast if $f$ is piecewise analytic.  Alternatively, if $f$ is smooth and the points $w$ are varied whilst the degrees $M$ are fixed, then the error decays like $h^{-M_{\min}-1}$, where $h = \max_{j=0,\ldots,k} | w_{j+1}  -w_j |$ and $M_{\min} = \min \{ M_0,\ldots,M_k \}$.

For analysis, we need to determine $\gamma_M$ and $\zeta_M$.  For the first we use the scaled Markov inequality $\| p' \|_{I} \leq \sqrt{2} M^2/ |I| \| p \|_{I}$, $\forall p \in \mathbb{P}_{M}$, $M \in \mathbb{N}$, where $|I|$ denotes the length of $I$.  Hence, if $\eta = \min_{j=0,\ldots,k} \{ w_{j+1} - w_j \}$ then $\gamma_M \leq \sqrt{2} M^2_{\max} / \eta$.  For $\zeta_M$, we recall the following inequality for polynomials $\| p \|_{\infty,I} \leq c M / \sqrt{|I|} \| p \|_{I}$, $\forall p \in \mathbb{P}_{M}$, $M \in \mathbb{N}$, where $c > 0$ is a constant.  Hence $\zeta_M \leq c M_{\max} / \sqrt{\eta}$.  We therefore deduce the following sufficient condition: $M^2_{\max} / \eta = \ord{K_N}$ as $N \rightarrow \infty$.  In the first scenario, where $\eta$ is fixed and $M_{\max}$ is varied, we attain the same square-root-type scaling for piecewise smooth functions when approximated by piecewise polynomials as with the polynomial space of the previous example.  In the second scenario, where $M_{\max}$ is fixed and $\eta$ is varied, we see that this leads to a linear relation between $K_N$ and $\eta$.  Thus, by forfeiting the superalgebraic/exponential convergence of the polynomial space for only algebraic convergence, we obtain a better scaling with $K_N$.  Note that in some cases it may be desirable to approximate using functions that are themselves smooth (up to a finite order).  In this case, we can replace $\rT_{w,M}$ by the spline space $\tilde{\rT}_{w,M_{\min}}$ of degree $M_{\min}$ on the knot sequence $w$.  Since $\tilde{\rT}_{w,M_{\min}} \subseteq \rT_{w,M}$ we obtain the same linear scaling with $K_N$ in this case as well.

\vspace{0.5pc} \noindent \textbf{Numerical results.} We demonstrate our results using two common nonuniform sampling schemes; jittered and log sampling (see \cite{BAMGACHNonuniform1D} for details).  In the first two panels of Fig.\ \ref{fig:rate_and_error_s}, we illustrate the scaling for different spaces $\rT_M$ between the sampling bandwidth $K_N$ and space dimension $M$ such that $C(N,M)$ is bounded.  For such $K_N$ and $M$, in the second pair of panels, we compute the $\rL^2$ error of the approximation $f_{N,M}$ for a continuous function $f$.  The superiority of the spline spaces for small $N$ is evident, with the polynomial space becoming better as $N$ increases.

\vspace{-0.5pc}
\begin{acknowledgement}
BA acknowledges support from the NSF DMS grant 1318894.  MG acknowledges support from the UK EPSRC grant EP/H023348/1 for the University of Cambridge Centre for Doctoral Training, the Cambridge Centre for Analysis.  AH acknowledges support from a Royal Society University Research Fellowship as well as the EPSRC grant EP/L003457/1.
\end{acknowledgement}

\vspace{-2pc}
\bibliographystyle{abbrv}
\small
\bibliography{1DRefs}

\begin{thebibliography}{10}

\bibitem{BAMGACHNonuniform1D}
B.~Adcock, M.~Gataric, and A.~C. Hansen.
\newblock On stable reconstructions from nonuniform {F}ourier measurements.
\newblock {\em SIAM J. Imaging Sci. (to appear)}, 2014.

\bibitem{BAACHAccRecov}
B.~Adcock and A.~C. Hansen.
\newblock {Stable reconstructions in {H}ilbert spaces and the resolution of the
  {G}ibbs phenomenon}.
\newblock {\em Appl. Comput. Harmon. Anal.}, 32(3):357--388, 2012.

\bibitem{AdcockHansenSpecData}
B.~Adcock and A.~C. Hansen.
\newblock Generalized sampling and the stable and accurate reconstruction of
  piecewise analytic functions from their {F}ourier coefficients.
\newblock {\em Math. Comp. (to appear)}, 2014.

\bibitem{BAACHOptimality}
B.~Adcock, A.~C. Hansen, and C.~Poon.
\newblock Beyond consistent reconstructions: optimality and sharp bounds for
  generalized sampling, and application to the uniform resampling problem.
\newblock {\em SIAM J. Math. Anal.}, 45(5):3114--3131, 2013.

\bibitem{AHPWavelet}
B.~Adcock, A.~C. Hansen, and C.~Poon.
\newblock On optimal wavelet reconstructions from {F}ourier samples: linearity
  and universality of the stable sampling rate.
\newblock {\em Appl. Comput. Harmon. Anal.}, 36(3):387--415, 2014.

\bibitem{AdcockHansenShadrinStabilityFourier}
B.~Adcock, A.~C. Hansen, and A.~Shadrin.
\newblock A stability barrier for reconstructions from {F}ourier samples.
\newblock {\em SIAM J. Numer. Anal.}, 52(1):125--139, 2014.

\bibitem{BoettcherMarkov}
A.~B{\"o}ttcher and P.~D{\"o}rfler.
\newblock Weighted {M}arkov-type inequalities, norms of {V}olterra operators,
  and zeros of {B}essel functions.
\newblock {\em Math. Nachr.}, 283(1):40--57, 2010.

\bibitem{GelbHinesEdge}
A.~Gelb and T.~Hines.
\newblock Detection of edges from nonuniform {F}ourier data.
\newblock {\em J. Fourier Anal. Appl. (to appear)}, 2011.

\bibitem{GelbHines2011Frames}
A.~Gelb and T.~Hines.
\newblock Recovering exponential accuracy from non-harmonic {F}ourier data
  through spectral reprojection.
\newblock {\em J. Sci. Comput.}, 51(158--182), 2012.

\bibitem{GelbSongFrameNFFT}
A.~Gelb and G.~Song.
\newblock {A frame theoretic approach to the Non-Uniform Fast Fourier
  Transform}.
\newblock {\em SIAM J. Numer. Anal. (to appear)}, 2014.

\bibitem{GrochenigIrregular}
K.~Gr{\"o}chenig.
\newblock {Reconstruction algorithms in irregular sampling}.
\newblock {\em Math. Comp.}, 59:181--194, 1992.

\bibitem{GrochenigIrregularExpType}
K.~Gr\"{o}chenig.
\newblock Irregular sampling, {T}oeplitz matrices, and the approximation of
  entire functions of exponential type.
\newblock {\em Math. Comp.}, 68(226):749--765, 1999.

\bibitem{hrycakIPRM}
T.~Hrycak and K.~Gr\"{o}chenig.
\newblock Pseudospectral {F}ourier reconstruction with the modified inverse
  polynomial reconstruction method.
\newblock {\em J. Comput. Phys.}, 229(3):933--946, 2010.

\bibitem{JacksonEtAlGridding}
J.~I. Jackson, C.~H. Meyer, D.~G. Nishimura, and A.~Macovski.
\newblock {Selection of a convolution function for {F}ourier inversion using
  gridding}.
\newblock {\em IEEE Trans. Med. Imaging}, 10:473--478, 1991.

\bibitem{MartinezEtAlEdge}
A.~Martinez, A.~Gelb, and A.~Gutierrez.
\newblock Edge detection from non-uniform {F}ourier data using the
  convolutional gridding algorithm.
\newblock {\em J. Sci. Comput. (to appear)}, 2014.

\bibitem{GelbPlatteEdge}
R.~Platte, A.~J. Gutierrez, and A.~Gelb.
\newblock Edge informed {F}ourier reconstruction from non-uniform spectral data
  with exponential convergence rates.
\newblock {\em Preprint}, 2012.

\bibitem{SzyldOblProj}
D.~Szyld.
\newblock The many proofs of an identity on the norm of oblique projections.
\newblock {\em Numer. Algorithms}, 42:309--323, 2006.

\bibitem{GelbNonuniformFourier}
A.~Viswanathan, A.~Gelb, D.~Cochran, and R.~Renaut.
\newblock {On reconstructions from non-uniform spectral data}.
\newblock {\em J. Sci. Comput.}, 45(1--3):487--513, 2010.

\end{thebibliography}

\end{document}